\documentclass[11pt]{article}
\usepackage[latin1]{inputenc}
\usepackage{amsmath,amsthm,amssymb}
\usepackage{amsfonts}
\usepackage{amsmath,amsthm,amssymb,amscd}
\usepackage{latexsym}
\usepackage{color}
\usepackage{graphicx}
\usepackage{mathrsfs}
\usepackage{cite}

\usepackage{color,enumitem,graphicx}
\usepackage[colorlinks=true,urlcolor=black,
citecolor=black,linkcolor=black,linktocpage,pdfpagelabels,
bookmarksnumbered,bookmarksopen]{hyperref}

\makeatletter \@addtoreset{equation}{section} \makeatother

\newtheorem{theorem}{Theorem}[section]

\newtheorem{proposition}{Proposition}[section]
\newtheorem{lemma}{Lemma}[section]
\newtheorem{remark}{Remark}[section]

\newtheorem{corollary}[theorem]{Corollary}

\allowdisplaybreaks

\begin{document}
\title{Normalized solutions for a biharmonic Choquard equation with exponential critical growth in $\mathbb{R}^4$}

\author{Wenjing Chen\footnote{Corresponding author.}\ \footnote{E-mail address:\, {\tt wjchen@swu.edu.cn} (W. Chen), {\tt zxwangmath@163.com} (Z. Wang).}\  \ and Zexi Wang\\
\footnotesize  School of Mathematics and Statistics, Southwest University,
Chongqing, 400715, P.R. China}

\date{ }
\maketitle

\begin{abstract}
{In this paper, we study the following biharmonic Choquard type problem
\begin{align*}
  \begin{split}
  \left\{
  \begin{array}{ll}
   \Delta^2u-\beta\Delta u=\lambda u+(I_\mu*F(u))f(u),
    \quad\mbox{in}\ \ \mathbb{R}^4,\\
    \displaystyle\int_{\mathbb{R}^4}|u|^2dx=c^2>0,\quad u\in H^2(\mathbb{R}^4),
    \end{array}
    \right.
  \end{split}
  \end{align*}
where $\beta\geq0$, $\lambda\in \mathbb{R}$, $I_\mu=\frac{1}{|x|^\mu}$ with $\mu\in (0,4)$,
$F(u)$ is the primitive function of $f(u)$, and $f$ is a continuous function with exponential critical growth.
We prove the existence of ground state normalized solutions for the above problem when the nonlinearity $f$ satisfies some conditions.}



\end{abstract}

\section{Introduction and statement of main result}
In this work, we are interested in the following biharmonic nonlinear Schr\"{o}dinger equation (NLS) with the mixed dispersion and a general Choquard nonlinear term 
\begin{equation}\label{back}
  i\partial_t\psi-\gamma\Delta^2\psi+\beta\Delta \psi+(I_\mu*F(\psi))f(\psi)=0, \quad\psi(0,x)=\psi_0(x),\quad \psi(t,x):\mathbb{R}\times \mathbb{R}^N\rightarrow \mathbb{C},
\end{equation}
where $N\geq1$, $i$ denotes the imaginary unit, $\gamma\geq0$, $\beta\in \mathbb{R}$, $I_\mu=\frac{1}{|x|^\mu}$ with $\mu\in [0,N)$, $\Delta^2:=\Delta(\Delta)$ is the biharmonic operator, and $f$ satisfies:


$(H_1)$ $ f(t)\in \mathbb{R}$ for $t\in \mathbb{R}$ and $f(e^{i\theta }z)=e^{i\theta }f(z)$ for any $\theta \in \mathbb{R}$, $z\in \mathbb{C}$;

$(H_2)$ $\displaystyle F(z)=\int_{0}^{|z|}f(t)dt$ for any $z\in \mathbb{C}$.

In physics, NLS is derived from the scalar nonlinear Helmhotz equation through the so-called paraxial approximation, see \cite{FIP}. Since its solutions may blow up at finite-time, it suggests that some terms neglected in the paraxial approximation may prevent this phenomenon. In \cite{FIP}, see also \cite{K,KS}, a small biharmonic term arises as a part of the nonparaxial correction to NLS, hence it is natural to consider this term as small but nonzero and study its effect on the blow up, which eventually gives rise to \eqref{back}.

If $\psi$ is a stationary wave, which means solution for \eqref{back} of the form $\psi(t,x)=e^{-i\lambda t}u(x)$ with $\lambda\in \mathbb{R}$ and $u\in H^2(\mathbb{R}^N)$ is a time-independent real value function, then $u$ satisfies
\begin{align}\label{e1.1}
   \gamma\Delta^2u-\beta\Delta u=\lambda u+(I_\mu*F(u))f(u),
    \quad\mbox{in}\ \ \mathbb{R}^N.
  \end{align}

A possible choice is to consider that $\lambda\in \mathbb{R}$ is given, the existence and qualitative properties of solutions for \eqref{e1.1} have been studied for a few decades by variational methods, see \cite{BS3,CGT,GS,MS1,MS2,MS3,MSV,RS1}.

From a physical point of view, it is interesting to find solutions of \eqref{e1.1} with a prescribed  $L^2$-norm,
and $\lambda\in \mathbb{R}$ appearing as a Lagrange multiplier. Solutions of this type are often referred to as normalized solutions. In fact, by $(H_1)$ and $(H_2)$,
we deduce that $F(z)\in \mathbb{R}$ and $f(z)\bar{z}\in \mathbb{R}$ for any $z\in \mathbb{C}$. Multiplying \eqref{back} by $\bar{\psi}$, integrating over $\mathbb{R}^N$ and taking the imaginary part, we find that $\frac{d}{dt}(\|\psi(t,\cdot)\|_2)=0$, this yields that $\|\psi(t,\cdot)\|_2=\|\psi(0,\cdot)\|_2$ for any $t\in \mathbb{R}$.

For a given $c>0$, let
 \begin{equation*}
  \widehat{S}(c):=\Big\{u\in H^2(\mathbb{R}^N):\int_{\mathbb{R}^N}|u|^2dx=c^2\Big\}.
\end{equation*}
  To obtain normalized solutions of \eqref{e1.1}, we find critical points of the functional
\begin{align*}
\widehat{\mathcal J}(u)
=\frac{\gamma}{2}\int_{\mathbb{R}^N}|\Delta u|^2dx+\frac{\beta}{2}\int_{\mathbb{R}^N}|\nabla u|^2dx-\frac{1}{2}\int_{\mathbb{R}^N}(I_\mu*F(u))F(u)dx
\end{align*}
on $\widehat{S}(c)$.
It is clear that any critical point $u_c$ of $\widehat{\mathcal J}$ on $\widehat{S}(c)$ corresponds to a Lagrange multiplier $\lambda_c\in \mathbb{R}$ such that $(u_c,\lambda_c)$ solves (weakly) \eqref{e1.1}. We call that $u$ is a ground state solution of \eqref{e1.1} on $\widehat{S}(c)$ if it minimizes the functional $\widehat{\mathcal J}$ on $\widehat{S}(c)$ among all the solutions.

When $\gamma=0$, $\beta=1$, and $\mu=0$ in \eqref{e1.1}, the problem reduces to the following nonlinear Schr\"{o}dinger equation with normalization constraint 
\begin{align}\label{mainuse}
  \begin{split}
  \left\{
  \begin{array}{ll}
   -\Delta u=\lambda u+f(u),
    \quad\mbox{in}\ \ \mathbb{R}^N,  \\
    \displaystyle\int_{\mathbb{R}^N}|u|^2dx=c^2.
    \end{array}
    \right.
  \end{split}
  \end{align}
If $f(u)=|u|^{p-2}u$, then $\bar{p}:=2+\frac{4}{N}$ is called the $L^2$-critical exponent, which comes from the Gagliardo-Nirenberg inequality \cite{Nirenberg1}.  For $p\in(2,2+\frac{4}{N})$ ($L^2$-subcritical case), the associated energy functional $\mathcal{J}_1$ of problem \eqref{mainuse} is bounded from below on
$
\widetilde{S}(c):=\{u\in H^1(\mathbb{R}^N):\int_{\mathbb{R}^N}|u|^2dx=c^2\}.
$
Thus, a ground state solution of \eqref{mainuse} can be found as a global minimizer of $\mathcal{J}_1$ on $\widetilde{S}(c)$.  Stuart \cite{S1,S2} first obtained the existence of normalized solutions of \eqref{mainuse} by the bifurcation theory. Other results for the $L^2$-subcritical problems can be found in \cite{CL,Lions2,Shibata} and references therein.

However, for $p\in[2+\frac{4}{N},2^*)$ ($L^2$-critical or supercritical case), where $2^*=\infty$ if $N\leq2$ and $2^*=\frac{2N}{N-2}$ if $N\geq3$, in this case, $\mathcal{J}_1$ is unbounded from below on $\widetilde{S}(c)$, and it seems impossible to search for a global minimizer of $\mathcal{J}_1$ on $\widetilde{S}(c)$. Furthermore, it is also difficult to get the boundedness and compactness of any $(PS)$ sequence. Jeanjean \cite{Jeanjean} first studied the existence of normalized solutions for \eqref{mainuse} in the $L^2$-supercritical case. Under the following conditions, Jeanjean considered an auxiliary functional $\widetilde{\mathcal{J}}_1:H^1(\mathbb{R}^N)\times \mathbb{R}\rightarrow \mathbb{R}$ defined by
 \begin{align*}
\widetilde{\mathcal{J}}_1(u,s)=\frac{e^{2s}}{2}\int_{\mathbb{R}^N}|\nabla u|^2dx-\frac{1}{e^{Ns}}\int_{\mathbb{R}^N}F(e^{\frac{Ns}{2}}u(x))dx.
\end{align*}

($W_1$) $f:\mathbb{R}\rightarrow\mathbb{R}$ is continuous and odd.

($W_2$) There exist $\alpha,\beta\in \mathbb{R}$ satisfying $2+\frac{4}{N}<\alpha\leq\beta<2^*$ such that
\begin{equation*}
  0<\alpha F(t)\leq f(t)t\leq\beta F(t)\quad\text{for any $t\neq0$, where $F(t)=\int_{0}^tf(s)ds$}.
\end{equation*}


($W_3$) Let $\widetilde{F}(t):=f(t)t-2F(t)$, then $\widetilde{F}'(t)$ exists and 
\begin{equation*}
  \widetilde{F}'(t)t>\Big(2+\frac{4}{N}\Big)\widetilde{F}(t)\quad\text{for any $t\neq0$}.
\end{equation*}
Under the conditions $(W_1)$ and $(W_2)$, Jeanjean proved that $\widetilde{\mathcal{J}}_1$ on $\widetilde{S}(c)\times \mathbb{R}$ has the same type of mountain-pass geometrical structure as $\mathcal{J}_1$ on $\widetilde{S}(c)$, and their mountain-pass levels $m(c)$ are equal. Then, using the Ekeland's variational principle and pseudo-gradient flow, the author obtained a bounded $(PS)_{m(c)}$ sequence $\{u_n\}$ for $\mathcal{J}_1$ on $\widetilde{S}(c)$. Next, by the Lagrange multiplier rule, the author proved that the Lagrange multiplier $\{\lambda_n\}$ is bounded in $\mathbb{R}$ and its limit $\lambda_c$ is strictly negative. Finally, it follows from the compact immersion of $H^1_{rad}(\mathbb{R}^N)\hookrightarrow L^p(\mathbb{R}^N)$ for any $p\in (2,2^*)$ that the nonlinearity is convergent, this fact together with $\lambda_c<0$ yields the compactness of $\{u_n\}$. In addition, it is worth mentioning that the author also proved the normalized solution is a ground state solution of \eqref{mainuse} if $(W_1)$, $(W_2)$ and $(W_3)$ hold.
The method of Jeanjean \cite{Jeanjean} has become a common method to deal with $L^2$-supercritical problems and has been developed by many researchers, see \cite{AJM,BM,BS0,BS2,Li1,DZ,Soave1,Soave2,WW,LZ2,JL} for normalized solutions in $\mathbb{R}^N$, \cite{NTV1,PV} for normalized solutions in bounded domains, and \cite{BLL,Li2,LY1,YCRS,YCT} for normalized solutions of Choquard equations. In particular, using the ideas of \cite{Jeanjean},
Alves et al. \cite{AJM} considered the existence of normalized solutions for \eqref{mainuse} in $\mathbb{R}^2$, where $f$ has an exponential critical growth.


When $\gamma=1$, $\beta\in \mathbb{R}$, $\mu=0$ in \eqref{e1.1} , the problem is related to the following mixed dispersion biharmonic NLS with prescribed $L^2$-norm constraint
\begin{align}\label{biharmonic}
  \begin{split}
  \left\{
  \begin{array}{ll}
  \Delta^2u-\beta\Delta u=\lambda u+f(u),
    \quad\mbox{in}\ \ \mathbb{R}^N,  \\
    \displaystyle\int_{\mathbb{R}^N}|u|^2dx=c^2,\quad u\in H^2(\mathbb{R}^N).
    \end{array}
    \right.
  \end{split}
  \end{align}
This kind of problem gives a new $L^2$-critical exponent $\bar{q}:=2+\frac{8}{N}$. 
 If $f(u)=|u|^{q-2}u$ and $q\in (2,2+\frac{8}{N})$, Bonheure et al. \cite{BCdN} studied \eqref{biharmonic} with $\beta>0$ and established the existence, qualitative properties of minimizers.
Using the profile decomposition of bounded sequences in $H^2(\mathbb{R}^N)$, Luo et al. \cite{LZZ} considered \eqref{biharmonic} with $\beta\in \mathbb{R}$ and obtained the existence, nonexistence, orbital stability of minimizers. Moreover,  Boussa\"{i}d et al. \cite{BFJ} improved the results of \cite{LZZ} by relaxing the extra restrictions on $\beta<0$ and $c$, which is totally new for $c>0$ small.

If $q\in [2+\frac{8}{N},4^*)$, it is impossible to search for a global minimizer to obtain a solution of \eqref{biharmonic}, where $4^*=\infty$ if $N\leq4$ and $4^*=\frac{2N}{N-4}$ if $N\geq5$. By using a minimax principle based on the homotopy stable family, Bonheure et al. \cite{BCGJ} obtained the existence of ground state solutions, radial positive solutions, and the multiplicity of radial solutions for \eqref{biharmonic} when $\beta>0$. Luo and Yang \cite{LY2} investigated the existence of at least two normalized solutions for \eqref{biharmonic} with $N\geq5$, $\beta<0$, where one is the local minimizer, the other one is the mountain-pass type solution.

Noting that all the aforementioned papers for \eqref{biharmonic} have only considered the homogeneous nonlinearity. With regard to this point, Luo and Zhang \cite{LZ1} first studied normalized solutions of \eqref{biharmonic} on $\widehat{S}(c)$ with a general nonlinear term $f$ and obtained the existence and orbital stability of minimizers,
where $\beta\in \mathbb{R}$ and $f$ satisfies the suitable $L^2$-subcritical assumptions.

Motivated by the results mentioned above, the aim of this paper is to study the existence of normalized solutions to the following biharmonic Choquard type equation
\begin{align}\label{e1.6}
  \begin{split}
  \left\{
  \begin{array}{ll}
   \Delta^2u-\beta\Delta u=\lambda u+(I_\mu*F(u))f(u),
    \quad\mbox{in}\ \ \mathbb{R}^4,  \\
    \displaystyle\int_{\mathbb{R}^4}|u|^2dx=c^2>0,\quad u\in H^2(\mathbb{R}^4),\\
    \end{array}
    \right.
  \end{split}
  \end{align}
where $\beta\geq0$, $\lambda\in \mathbb{R}$, $I_\mu=\frac{1}{|x|^\mu}$ with $\mu\in(0,4)$,
and $f$ has exponential critical growth in the sense of the Adams inequality \cite{RS}.
We assume that $f$ satisfies:

$(f_1)$ $f:\mathbb{R}\rightarrow \mathbb{R}$ is continuous and odd, $\displaystyle\lim\limits_{t\rightarrow0} \frac{|f(t)|}{|t|^\nu}=0$ for some $\nu>2-\frac{\mu}{4}$;

$(f_2)$ $f$ has exponential critical growth at infinity. Precisely, it holds
\begin{align*}
  \begin{split}
  \lim\limits_{|t|\rightarrow+\infty}\frac{|f(t)|}{e^{\alpha t^2}}=\left\{
  \begin{array}{ll}
  0,&\ \ \text {for} \,\,\,\alpha>32\pi^2,\\
  +\infty,&\ \  \text {for} \,\,\,0<\alpha<32\pi^2;
    \end{array}
    \right.
  \end{split}
  \end{align*}

$(f_3)$ There exists a constant $\theta>3-\frac{\mu}{4}$ such that
\begin{equation*}
  0<\theta F(t)\leq tf(t), \,\,\,\text{for all $t\in \mathbb{R}\backslash \{0\}$};
\end{equation*}

$(f_4)$ There exist constants $p>3-\frac{\mu}{4}$ and $\tau>0$ such that
\begin{equation*}
  F(t)\geq \frac{\tau}{p} |t|^p,\,\,\,\text{for all $t\in \mathbb{R}$};
\end{equation*}

$(f_5)$ For any $t\in \mathbb{R}\backslash \{0\}$, let $\overline{F}(t):=f(t)t-(2-\frac{\mu}{4})F(t)$, $f'(t)$
exists, and
\begin{equation*}
(3-\frac{\mu}{4})F(s)\overline{F}(t)<F(s)\overline{F}'(t)t+\overline{F}(s)(\overline{F}(t)-F(t)),
 \quad \text{for any $s,t\in \mathbb{R}\backslash \{0\}$};
\end{equation*}

$(f_6)$ $\frac{\overline{F}(t)}{|t|^{3-\frac{\mu}{4}}}$ is non-increasing in $(-\infty,0)$ and non-decreasing in $(0,+\infty)$.

Our main result is as follows:
\begin{theorem}\label{th1}
Assume that $c\in (0,\sqrt{\frac{4-\mu}{4}})$ and $f$ satisfies $(f_1)-(f_4)$,
$(f_5)$ or $(f_6)$, then there exists $\tau^*>0$ such that for any $\tau\geq\tau^*$, problem \eqref{e1.6} possesses a non-negative radial ground state solution.
\end{theorem}

\begin{remark}
{\rm
A typical example satisfying $(f_1)-(f_6)$ is
\begin{equation*}
  f(t)=\tau|t|^{p-2}te^{32\pi t^2}
\end{equation*}
for any $p>\max\{3, \nu+1,\theta \}$. If $f'(t)$ exists, the condition $(f_6)$ implies $(f_5)$. Indeed, since
\begin{equation*}
  \frac{d}{dt}\Big(\frac{\overline{F}(t)}{|t|^{3-\frac{\mu}{4}}}\Big)=\frac{\overline{F}'(t)t-(3-\frac{\mu}{4})\overline{F}(t)}{|t|^{4-\frac{\mu}{4}}sign(t)},
\end{equation*}
and thus by $(f_6)$, we have $\overline{F}'(t)t-(3-\frac{\mu}{4})\overline{F}(t)\geq 0$ for any $t\in \mathbb{R}\backslash \{0\}$. From $(f_3)$, we know $F(t)>0$, $\overline{F}(t)>0$, and $\overline{F}(t)-F(t)>0$ for any $t\in \mathbb{R}\backslash \{0\}$, hence the condition $(f_5)$ follows.
}
\end{remark}
Let us define the energy functional $\mathcal J:H^2(\mathbb{R}^4)\rightarrow\mathbb{R}$ by
\begin{align*}
\mathcal J(u)=\frac{1}{2}\int_{\mathbb{R}^4}|\Delta u|^2dx+\frac{\beta}{2}\int_{\mathbb{R}^4}|\nabla u|^2dx-\frac{1}{2}\int_{\mathbb{R}^4}(I_\mu*F(u))F(u)dx,
\end{align*}
where $H^2(\mathbb{R}^4)=\Big\{u\in L^2(\mathbb{R}^4):\nabla u\in L^2(\mathbb{R}^4),\ \Delta u\in L^2(\mathbb{R}^4)\Big\}$.

By $(f_1)$ and $(f_2)$, fix $q>0$, for any $\xi>0$ and $\alpha>32\pi^2$, there exists a constant $C_\xi>0$ such that
\begin{equation*}\label{fcondition1}
  |f(t)|\leq \xi|t|^\nu+C_\xi|t|^q(e^{\alpha t^2}-1)\quad\text{for all $t\in \mathbb{R}$},
\end{equation*}
and using $(f_3)$, we have
\begin{equation}\label{fcondition2}
|F(t)|\leq \xi|t|^{\nu+1}+C_\xi|t|^{q+1}(e^{\alpha t^2}-1)\quad\text{for all $t\in \mathbb{R}$}.
\end{equation}
By \eqref{fcondition2}, using the Hardy-Littlewood-Sobolev inequality \cite{LL} and the Adams inequality \cite{RS,Y}, we obtain
$\mathcal{J}$ is well defined in $H^2(\mathbb{R}^4)$ and $\mathcal{J}\in C^1(H^2(\mathbb{R}^4),\mathbb{R})$ with
\begin{align*}
\langle\mathcal J'(u),v\rangle=&\int_{\mathbb{R}^4}\Delta u \Delta vdx+\beta\int_{\mathbb{R}^4}\nabla u \cdot \nabla vdx-\int_{\mathbb{R}^4}(I_\mu*F(u))f(u)vdx\\
=&\int_{\mathbb{R}^4}v\Delta^2 u dx-\beta\int_{\mathbb{R}^4}v\Delta u dx-\int_{\mathbb{R}^4}(I_\mu*F(u))f(u)vdx,
\end{align*}
for any $u, v\in H^2(\mathbb{R}^4)$.

For any $c\in (0,\sqrt{\frac{4-\mu}{4}})$, set
\begin{equation*}
 S(c):=\Big\{u\in H^2(\mathbb{R}^4):\int_{\mathbb{R}^4}|u|^2dx=c^2\Big\}.
\end{equation*}
Define
\begin{equation*}
  E(c):=\inf\limits_{u\in \mathcal{P}(c)}\mathcal J(u),
\end{equation*}
where $\mathcal{P}(c)$ is the Pohozaev manifold defined by $\mathcal{P}(c)=\Big\{u\in S(c):P(u)=0\Big\}$ with
\begin{align*}
  P(u)&=2\int_{\mathbb{R}^4}|\Delta u |^2\ dx+\beta\int_{\mathbb{R}^4}|\nabla u|^2 dx+\frac{8-\mu}{2} \int_{\mathbb{R}^4} (I_\mu*F(u))F(u)dx-2\int_{\mathbb{R}^4} (I_\mu*F(u))f(u)u dx.
\end{align*}
As will be shown in Lemma \ref{equi}, $\mathcal{P}(c)$ is nonempty, and from Lemma \ref{Pohozaev}, we can see that any critical point of $\mathcal{J}$ on $S(c)$ stays in $\mathcal{P}(c)$, thus any critical point $u$ of $\mathcal{J}$ on $S(c)$ with $\mathcal J(u)=E(c)$ is a ground state solution of \eqref{e1.6}.

For any $s\in \mathbb{R}$ and $u\in H^2(\mathbb{R}^4)$, we define
\begin{equation*}
  \mathcal{H}(u,s)(x)=e^{2s}u(e^sx),\quad \text{for a.e. $x\in \mathbb{R}^4$}.
\end{equation*}
For simplicity, we always write $\mathcal{H}(u,s)$. One can easily check that $\| \mathcal{H}(u,s)\|_2=\|u\|_2$ for any $s\in \mathbb{R}$.

\begin{remark}
{\rm
$(i)$ Assumption $(f_3)$ is the classical Ambrosetti-Rabinowitz condition, which is used to guarantee the boundedness of $(PS)_{m_\tau(c)}$ sequence. Moreover, we need $(f_3)$ to establish the relationship between the mountain-pass level $m_\tau(c)$ and the ground state energy $E(c)$.

$(ii)$ $(f_4)$ is crucial in our approach to obtain the mountain-pass geometrical structure and an upper bound of $m_\tau(c)$ when $\tau>0$ large enough.

$(iii)$ Under the condition $(f_5)$ or $(f_6)$, we can prove the uniqueness of $s_u\in \mathbb{R}$ such that $\mathcal{H}(u,s_u)\in \mathcal{P}(c)$, for any $u\in S(c)$, this together with $(f_3)$ yields $m_\tau(c)=E(c)$.}
\end{remark}

\begin{remark}
{\rm To prove Theorem \ref{th1}, we use the ideas developed in \cite{Jeanjean} to construct a bounded $(PS)_{m_\tau(c)}$ sequence, then proving the compactness of the $(PS)_{m_\tau(c)}$ sequence, we obtain the existence of normalized solutions of problem \eqref{e1.6}.  The main difficulty in the proof is the lack of compactness for $(PS)_{m_\tau(c)}$ sequence due to the whole space and exponential critical growth. In order to overcome this difficulty, we work directly in the space of radially symmetric functions in $H^2(\mathbb{R}^4)$ and use $(f_4)$, the Adams inequality. It is worth pointing out that the key step to use the Adams inequality is to give a suitable uniformly upper bound on the $H^2$
-norm for the $(PS)_{m_\tau(c)}$ sequence. 
 Finally, by $(f_3)$, $(f_5)$ or $(f_3)$, $(f_6)$, we can prove that the normalized solution is a ground state solution of \eqref{e1.6}.}
\end{remark}
\begin{remark}
{\rm Compared with \cite{AJM} and \cite{Jeanjean}, 
the appearance of the nonlocal convolution term brings new difficulties when we prove the convergence of the nonlinearity, see Lemmas \ref{strong2} and \ref{strong3}. Using the
Hardy-Littlewood-Sobolev inequality and the Adams inequality, we prove that the nonlocal convolution term is bounded a.e. in $\mathbb{R}^4$. This fact transforms nonlocal problem into a local one, then using a variant of the Lebesgue dominated convergence theorem, we overcome this difficulty. For this reason, we shall restrict $c\in (0,\sqrt{\frac{4-\mu}{4}})$.}
\end{remark}

We will often make use of the following interpolation inequality
\begin{equation}\label{ineq}
  \int_{\mathbb{R}^4}|\nabla u|^2dx\leq  \Big(\int_{\mathbb{R}^4}|\Delta u|^2dx\Big)^{\frac{1}{2}}\Big(\int_{\mathbb{R}^4}| u|^2dx\Big)^{\frac{1}{2}},
\end{equation}
and the well-known Gagliardo-Nirenberg inequality \cite{Nirenberg1}, namely, for $u\in H^2(\mathbb{R}^4)$ and $p\geq2$,
\begin{equation}\label{GNinequality}
  \|u\|_p\leq B_{p}\|\Delta u\|_2^{\frac{p-2}{p}}\| u\|_2^{\frac{2}{p}},
\end{equation}
where $B_{p}$ is a constant depending on $p$.

In this paper, $C$ denotes positive constant possibly different from line to line. $L^p(\mathbb{R}^4)$ is the usual Lebesgue space endowed with the norm $\|u\|_{p}=(\int_{\mathbb{R}^4}|u|^pdx)^{\frac{1}{p}}$ when $1<p<\infty$, $\|u\|_\infty=\inf\{C>0, |u(x)|\leq C \,\,\text{a.e. in} \,\,\mathbb{R}^4\}$. $H^2_{rad}(\mathbb{R}^4)$ is the space of radially symmetric functions in $H^2(\mathbb{R}^4)$ and its corresponding norm is defined by
$$\|u\|:=(\|\Delta u\|_2^2+2\|\nabla u\|_2^2+\|u\|_2^2)^{\frac{1}{2}}.
$$
From \eqref{ineq}, we can see that $\|\cdot\|$ is equivalent to the norm $\|u\|_{H^2(\mathbb{R}^4)}:=(\|\Delta u\|_2^2+\|u\|_2^2)^{\frac{1}{2}}$.

The paper is organized as follows. Section \ref{sec preliminaries} contains some useful lemmas.
In Section \ref{minimax}, we use the mountain-pass argument to construct a bounded $(PS)_{m_\tau(c)}$ sequence.
Section \ref{main} is devoted to the proof of Theorem \ref{th1}.

\section{Preliminaries}\label{sec preliminaries}

In this section, we give some preliminaries.
For the nonlocal type problems with Riesz potential, an important inequality due to the Hardy-Littlewood-Sobolev inequality will be used in the following.
\begin{proposition}\cite[Theorem 4.3]{LL}\label{HLS}
Assume that $1<r$, $t<\infty$, $0<\mu<4$ and
$
\frac{1}{r}+\frac{\mu}{4}+\frac{1}{t}=2.
$
Then there exists $C(\mu,r,t)>0$ such that
\begin{align}\label{HLSin}
\Big|\int_{\mathbb R^{4}}(I_{\mu}\ast g(x))h(x)dx\Big|\leq
C(\mu,r,t)\|g\|_r
\|h\|_t,
\end{align}
for all $g\in L^r(\mathbb{R}^4)$ and $h\in L^t(\mathbb{R}^4)$.
If $ r=t=\frac{8}{8-\mu}$, equality holds in \eqref{HLSin} if and only if $h=cg$ for a constant $c$ and
\begin{equation*}
  g(x)=A\Big(\frac{1}{m^2+|x-a|^2}\Big)^{\frac{8-\mu}{2}}
\end{equation*}
for some $A\in \mathbb{C}$, $0\neq m\in \mathbb{R}$ and $a\in \mathbb{R}^4$.
\end{proposition}
From Proposition \ref{HLS}, we have the following result.
\begin{corollary}\label{remarkHLS}
Let $\mu\in(0,4)$ and $s\in(1,\frac{4}{4-\mu})$. If $\varphi\in L^s(\mathbb{R}^4)$, then $I_\mu\ast\varphi\in L^{\frac{4s}{4-(4-\mu)s}}(\mathbb{R}^4)$ and
\begin{equation*}
  \|I_\mu\ast\varphi\|_{\frac{4s}{4-(4-\mu) s}}\leq C(\mu, s)
  \|\varphi\|_s ,
\end{equation*}
where the constant $C( \mu, s)>0$ depends on $\mu, s$.
\end{corollary}
\begin{proof}
Denote $T\varphi:=I_\mu\ast\varphi$ and $t:=\frac{4s}{(8-\mu)s-4}$, then $\frac{1}{t}+\frac{4-(4-\mu) s}{4s}=1$ and
\begin{equation*}
  \|T\varphi\|_{\frac{4s}{4-(4-\mu) s}}=\sup\limits_{\|h\|_t=1}|\langle T\varphi,h\rangle|=\sup\limits_{\|h\|_t=1}\Big|\int_{\mathbb R^{4}}(I_{\mu}\ast \varphi(x))h(x)dx\Big|.
\end{equation*}
Since $\varphi\in L^s(\mathbb{R}^4)$, $h\in L^t(\mathbb{R}^4)$ and $\frac{1}{s}+\frac{\mu}{4}+\frac{1}{t}=2$, by \eqref{HLSin}, we deduce that
\begin{equation*}
  \|T\varphi\|_{\frac{4s}{4-(4-\mu) s}}\leq C(\mu,s)\|\varphi\|_s.
\end{equation*}
\end{proof}

\begin{lemma}\cite[Lemma 4.8]{Kav}\label{weakcon}
Let $\Omega\subseteq\mathbb R^{4}$ be any open set. For $1<s<\infty$, let $\{u_n\}$ be bounded in $L^s(\Omega)$ and $u_n(x)\rightarrow u(x)$ a.e. in $\Omega$. Then $u_n(x)\rightharpoonup u(x) \,\,in \,\,L^s(\Omega)$.
\end{lemma}

We recall the following Adams type inequality.
\begin{lemma}\cite{RS,Y}\label{adams}
(i) If $\alpha>0$ and $u\in H^2(\mathbb{R}^4)$, then
\begin{equation*}
  \int_{\mathbb{R}^4}(e^{\alpha u^2}-1)dx<+\infty;
\end{equation*}
(ii) There exists a constant $C>0$ such that
\begin{equation*}
  \sup\limits_{u\in H^2(\mathbb{R}^4), \|u\|\leq1}\int_{\mathbb{R}^4}(e^{\alpha u^2}-1)dx\leq C
\end{equation*}
for all $\alpha\leq 32\pi^2$.
This inequality is sharp, in the sense that if $\alpha>32\pi^2$, then the supremum is infinite.
\end{lemma}

\begin{lemma}\label{unibdd}
Suppose that $c\in (0,\sqrt{\frac{4-\mu}{4}})$, let $\{u_n\}\subset S(c)$ be a sequence 
satisfying $\limsup\limits_{n\rightarrow\infty}\|\Delta u_n\|_2<\sqrt{\frac{4-\mu}{4}}-c$. Then up to a subsequence, there exist  $\alpha>32\pi^2$ close to $32\pi^2$, $t>1$ close to $1$ and $C>0$ such that
\begin{equation*}
 \sup\limits_{n\in \mathbb{N}^+} \int_{\mathbb{R}^4}(e^{\alpha u_n^2}-1)^tdx\leq C.
\end{equation*}
\end{lemma}
\begin{proof}
By $\limsup\limits_{n\rightarrow\infty}\|\Delta u_n\|_2<\sqrt{\frac{4-\mu}{4}}-c$, up to a subsequence, we assume that $\sup\limits_{n\in \mathbb{N}^+}\|\Delta u_n\|_2<\sqrt{\frac{4-\mu}{4}}-c$. From \eqref{ineq} and $\{u_n\} \subset S(c)$, we have  $$\sup\limits_{n\in \mathbb{N}^+}\|u_n\|^2\leq \sup\limits_{n\in \mathbb{N}^+}\|\Delta u_n\|_2^2+2c\sup\limits_{n\in \mathbb{N}^+}\|\Delta u_n\|_2+c^2<\frac{4-\mu}{4}<1,$$ thus there exists $m\in(0,1)$ such that
\begin{equation*}
 \sup\limits_{n\in \mathbb{N}^+} \|u_n\|^2<m<1.
\end{equation*}
Fix $\alpha>32\pi^2$ close to $32\pi^2$ and $t>1$ close to $1$ such that $\alpha tm\leq 32\pi^2$, it yields that
\begin{equation*}
 \sup\limits_{n\in \mathbb{N}^+} \int_{\mathbb{R}^4}(e^{\alpha u_n^2}-1)^tdx\leq \sup\limits_{n\in \mathbb{N}^+}\int_{\mathbb{R}^4}(e^{\alpha t m(\frac{u_n}{\|u_n\|})^2}-1)dx.
\end{equation*}
Hence, by Lemma \ref{adams}, there exists $C>0$ such that
\begin{equation*}
  \sup\limits_{n\in \mathbb{N}^+}\int_{\mathbb{R}^4}(e^{\alpha u_n^2}-1)^tdx\leq C.
\end{equation*}
\end{proof}
\begin{lemma}\label{strong1}
Suppose that $c\in (0,\sqrt{\frac{4-\mu}{4}})$, let $\{u_n\}\subset S_r(c):=S(c)\cap H^2_{rad}(\mathbb{R}^4)$ be a sequence 
 and satisfy $\limsup\limits_{n\rightarrow\infty}\|\Delta u_n\|_2<\sqrt{\frac{4-\mu}{4}}-c$. Up to a subsequence, if $u_n\rightharpoonup u$ in $H^2_{rad}(\mathbb{R}^4)$, then there exists $\alpha>32\pi^2$ close to $32\pi^2$ such that for all $q>0$,
\begin{equation*}
  \int_{\mathbb{R}^4}|u_n|^{q+1}(e^{\alpha u_n^2}-1)dx\rightarrow \int_{\mathbb{R}^4}|u|^{q+1}(e^{\alpha u^2}-1)dx,\,\,\,\text{as}\,\,\, n\rightarrow\infty.
\end{equation*}
\end{lemma}

\begin{proof}
By Lemma \ref{unibdd}, there exist $\alpha>32\pi^2$ close to $32\pi^2$, $t>1$ close to $1$ such that $(e^{\alpha u_n^2}-1)$ is uniformly bounded in $L^t(\mathbb{R}^4)$. Since $u_n\rightharpoonup u$ in $H^2_{rad}(\mathbb{R}^4)$, then $u_n\rightarrow u$ a.e. in $\mathbb{R}^4$. By Lemma \ref{weakcon}, we obtain $(e^{\alpha u_n^2}-1)\rightharpoonup (e^{\alpha u^2}-1)$ in $L^t(\mathbb{R}^4)$. Moreover, for $t'=\frac{t}{t-1}$, using the compact  embedding $H^2_{rad}(\mathbb{R}^4)\hookrightarrow L^{(q+1)t'}(\mathbb{R}^4)$, we derive that $u_n\rightarrow u$ in $L^{(q+1)t'}(\mathbb{R}^4)$ as $n\rightarrow\infty$, and so
 $|u_n|^{q+1}\rightarrow |u|^{q+1}$ in $L^{t'}(\mathbb{R}^4)$ as $n\rightarrow\infty$, i.e.,
 \begin{equation*}
  \Big(\int_{\mathbb{R}^4}\Big||u_n|^{q+1}-|u|^{q+1}\Big|^{t'}dx\Big)^{\frac{1}{t'}}\rightarrow 0,\,\,\,\text{as}\,\,\, n\rightarrow\infty.
\end{equation*}
Thus, by the definition of weak convergence, using the H\"{o}lder inequality, we infer that
\begin{align*}
  &\Big|\int_{\mathbb{R}^4}|u_n|^{q+1}(e^{\alpha u_n^2}-1)dx-\int_{\mathbb{R}^4}|u|^{q+1}(e^{\alpha u^2}-1)dx\Big|\\
&\leq\int_{\mathbb{R}^4}\Big||u_n|^{q+1}-|u|^{q+1}\Big|(e^{\alpha u_n^2}-1)dx+\int_{\mathbb{R}^4}|u|^{q+1}\Big|(e^{\alpha u_n^2}-1)-(e^{\alpha u^2}-1)\Big|dx\\
&\leq  \Big(\int_{\mathbb{R}^4}\Big||u_n|^{q+1}-|u|^{q+1}\Big|^{t'}dx\Big)^{\frac{1}{t'}}\Big(\int_{\mathbb{R}^4}(e^{\alpha u_n^2}-1)^tdx\Big)^{\frac{1}{t}}\\
  &\quad+\int_{\mathbb{R}^4}|u|^{q+1}\Big|(e^{\alpha u_n^2}-1)-(e^{\alpha u^2}-1)\Big|dx\rightarrow 0,
\end{align*}
as $n\rightarrow\infty$.
\end{proof}

\begin{lemma}\label{strong2}
Assume that $(f_1)$, $(f_2)$ hold and $c\in (0,\sqrt{\frac{4-\mu}{4}})$. Let $\{u_n\}\subset S_r(c)$ be a sequence 
and satisfy $\limsup\limits_{n\rightarrow\infty}\|\Delta u_n\|_2<\sqrt{\frac{4-\mu}{4}}-c$. Up to a subsequence, if $u_n\rightharpoonup u$ in $H^2_{rad}(\mathbb{R}^4)$, then 
\begin{equation*}
  \int_{\mathbb{R}^4}(I_\mu*F(u_n))f(u_n)u_ndx\rightarrow \int_{\mathbb{R}^4}(I_\mu*F(u))f(u)udx,\,\,\,\text{as}\,\,\, n\rightarrow\infty.
\end{equation*}
\end{lemma}

\begin{proof}
By Corollary \ref{remarkHLS}, let $s\rightarrow \frac{4}{4-\mu}$, then $\frac{4s}{4-(4-\mu)s}\rightarrow\infty$, thus we have
\begin{equation*}
  \|I_\mu\ast F(u_n)\|_\infty\leq C\|F(u_n)\|_{\frac{4}{4-\mu}}\quad\text{for all $n\in \mathbb{N}^+$}.
\end{equation*}
and by \eqref{fcondition2}, fix $q>0$, for any $\xi>0$ and $\alpha>32\pi^2$ close to $32\pi^2$, there exists a constant $C_\xi>0$
such that
\begin{equation*}
  |F(u_n)|\leq \xi|u_n|^{\nu+1}+C_\xi|u_n|^{q+1}(e^{\alpha u_n^2}-1)\quad\text{for all $n\in \mathbb{N}^+$}.
\end{equation*}
By Lemma \ref{unibdd}, up to a subsequence, there exists $m\in (0,\frac{4-\mu}{4})$ such that $\sup\limits_{n\in \mathbb{N}^+}\|u_n\|^2<m$.
Fix $t>1$ close to $1$ such that $\frac{4\alpha mt}{4-\mu}\leq 32\pi^2$, for $t'=\frac{t}{t-1}$, using the H\"{o}lder inequality, the Sobolev inequality and Lemma \ref{adams}, we have
\begin{align*}
  \|F(u_n)\|_{\frac{4}{4-\mu}}&\leq\bigg(\int_{\mathbb{R}^4}\Big(\xi|u_n|^{\nu+1}+C_\xi|u_n|^{q+1}(e^{\alpha u_n^2}-1)\Big)^{\frac{4}{4-\mu}}dx\bigg)^{\frac{4-\mu}{4}}\nonumber\\
  &\leq\bigg(\int_{\mathbb{R}^4}\Big(C|u_n|^{\frac{4(\nu+1)}{4-\mu}}+C|u_n|^{\frac{4(q+1)}{4-\mu}}(e^{\frac{4\alpha u_n^2}{4-\mu}}-1)\Big)dx\bigg)^{\frac{4-\mu}{4}}\nonumber\\
  &\leq C\|u_n\|_{\frac{4(\nu+1)}{4-\mu}}^{\nu+1}+C\bigg(\int_{\mathbb{R}^4}|u_n|^{\frac{4(q+1)}{4-\mu}}(e^{\frac{4\alpha u_n^2}{4-\mu}}-1)dx\bigg)^{\frac{4-\mu}{4}}\nonumber\\
  &\leq C\|u_n\|_{\frac{4(\nu+1)}{4-\mu}}^{\nu+1}+C\Big(\int_{\mathbb{R}^4}|u_n|^{\frac{4(q+1)t'}{4-\mu}}dx\Big)^{\frac{4-\mu}{4t'}}\Big(\int_{\mathbb{R}^4}(e^{\frac{4\alpha mt}{4-\mu}(\frac{u_n}{\|u_n\|})^2}-1)dx\Big)^{\frac{4-\mu}{4t}}\nonumber\\
  &\leq C\|u_n\|_{\frac{4(\nu+1)}{4-\mu}}^{\nu+1}+C\|u_n\|_{\frac{4(q+1)t'}{4-\mu}}^{q+1}\leq C\|u_n\|^{\nu+1}+C\|u_n\|^{q+1}\leq C.
\end{align*}
Hence
\begin{equation*}
  (I_\mu*F(u_n))f(u_n)u_n\rightarrow (I_\mu*F(u))f(u)u\quad \text{a.e. in $\mathbb{R}^4$},
\end{equation*}
\begin{equation*}
  |(I_\mu*F(u_n))f(u_n)u_n|\leq C|f(u_n)u_n|\leq C|u_n|^{\nu+1}+C|u_n|^{q+1}(e^{\alpha u_n^2}-1)\quad\text{for all $n\in \mathbb{N}^+$},
\end{equation*}
and
\begin{equation*}
  |u_n|^{\nu+1}+|u_n|^{q+1}(e^{\alpha u_n^2}-1)\rightarrow |u|^{\nu+1}+|u|^{q+1}(e^{\alpha u^2}-1)\quad\text{a.e. in $\mathbb{R}^4$}.
\end{equation*}
By Lemma \ref{strong1}, we obtain
\begin{equation*}
  \int_{\mathbb{R}^4}|u_n|^{q+1}(e^{\alpha u_n^2}-1)dx\rightarrow \int_{\mathbb{R}^4}|u|^{q+1}(e^{\alpha u^2}-1)dx,\,\,\,\text{as}\,\,\, n\rightarrow\infty.
\end{equation*}
By the compact embedding $H^2_{rad}(\mathbb{R}^4)\hookrightarrow L^{\nu+1}(\mathbb{R}^4)$, one has
\begin{equation*}
  \int_{\mathbb{R}^4}|u_n|^{\nu+1}dx\rightarrow \int_{\mathbb{R}^4}|u|^{\nu+1}dx,\,\,\,\text{as}\,\,\, n\rightarrow\infty.
\end{equation*}
Applying a variant of the Lebesgue dominated convergence theorem, we get
\begin{equation*}
  \int_{\mathbb{R}^4}(I_\mu*F(u_n))f(u_n)u_ndx\rightarrow \int_{\mathbb{R}^4}(I_\mu*F(u))f(u)udx,\,\,\,\text{as}\,\,\, n\rightarrow\infty.
\end{equation*}

\end{proof}
\begin{lemma}\label{strong3}
Assume that $(f_1)$, $(f_2)$ hold and $c\in (0,\sqrt{\frac{4-\mu}{4}})$. Let $\{u_n\}\subset S(c)$ be a sequence 
satisfying $\limsup\limits_{n\rightarrow\infty}\|\Delta u_n\|_2<\sqrt{\frac{4-\mu}{4}}-c$. Then up to a subsequence, if $u_n\rightharpoonup u$ in $H^2(\mathbb{R}^4)$, for any $\phi\in C_0^\infty(\mathbb{R}^4)$, we have
\begin{equation}\label{cov1}
  \int_{\mathbb{R}^4}\Delta u_n\Delta\phi dx \rightarrow  \int_{\mathbb{R}^4}\Delta u\Delta\phi dx,\,\,\,\text{as}\,\,\, n\rightarrow\infty,
\end{equation}
\begin{equation}\label{cov2}
  \int_{\mathbb{R}^4}\nabla u_n\cdot\nabla\phi dx \rightarrow  \int_{\mathbb{R}^4}\nabla u\cdot\nabla\phi dx,\,\,\,\text{as}\,\,\, n\rightarrow\infty,
\end{equation}
\begin{equation}\label{cov3}
  \int_{\mathbb{R}^4}u_n\phi dx \rightarrow  \int_{\mathbb{R}^4}u\phi dx,\,\,\,\text{as}\,\,\, n\rightarrow\infty,
\end{equation}
and
\begin{equation}\label{cov4}
  \int_{\mathbb{R}^4} (I_\mu*F(u_n))f(u_n)\phi dx\rightarrow \int_{\mathbb{R}^4} (I_\mu*F(u))f(u)\phi dx,\,\,\,\text{as}\,\,\, n\rightarrow\infty.
\end{equation}
\end{lemma}

\begin{proof}
For any fixed $v\in H^2(\mathbb{R}^4)$, define
\begin{equation*}
  f_v(u):=\int_{\mathbb{R}^4}\Delta u\Delta v dx,\quad g_v(u):=\int_{\mathbb{R}^4}\nabla u\cdot\nabla v dx,\quad h_v(u):=\int_{\mathbb{R}^4}u v dx\quad \text{for every $u \in H^2(\mathbb{R}^4)$}.
\end{equation*}
By the H\"{o}lder inequality, we have
\begin{equation*}
  |f_v(u)|,|g_v(u)|,|h_v(u)|\leq\|v\|\|u\|.
\end{equation*}
This yields that $f_v$, $g_v$, and $h_v$ are continuous linear functionals on $H^2(\mathbb{R}^4)$. Thus, by $u_n\rightharpoonup u$ in $H^2(\mathbb{R}^4)$ and $C_0^\infty(\mathbb{R}^4)$ is dense in $H^2(\mathbb{R}^4)$, we obtain that (\ref{cov1})-(\ref{cov3}) hold.

Next, we prove (\ref{cov4}). By Lemma \ref{strong1},
we know that there exist $\alpha>32\pi^2$ close to $32\pi^2$, $t>1$ close to $1$ such that $(e^{\alpha u_n^2}-1)$ is uniformly bounded in $L^t(\mathbb{R}^4)$, and $(e^{\alpha u_n^2}-1)\rightharpoonup (e^{\alpha u^2}-1)$ in $L^t(\mathbb{R}^4)$.
By Lemma \ref{strong2}, there exists a constant $C>0$ such that  $\|I_\mu\ast F(u_n)\|_\infty\leq C$ for all $n\in \mathbb{N}^+$.
Hence, for any $\phi\in C_0^\infty(\mathbb{R}^4)$, we have
\begin{equation*}
  (I_\mu*F(u_n))f(u_n)\phi\rightarrow (I_\mu*F(u))f(u)\phi\quad \text{a.e. in $\mathbb{R}^4$},
\end{equation*}
\begin{equation*}
  |(I_\mu*F(u_n))f(u_n)\phi|\leq C|f(u_n)||\phi|\leq C|u_n|^{\nu}|\phi|+C|u_n|^{q}|\phi|(e^{\alpha u_n^2}-1)\quad\text{for all $n\in \mathbb{N}^+$},
\end{equation*}
and
\begin{equation*}
  |u_n|^{\nu}|\phi|+|u_n|^{q}|\phi|(e^{\alpha u_n^2}-1)\rightarrow |u|^{\nu}|\phi|+|u|^{q}|\phi|(e^{\alpha u^2}-1)\quad\text{a.e. in $\mathbb{R}^4$}.
\end{equation*}
Denote $\Omega=supp \phi$. In the following, we prove that
 \begin{equation}\label{strong3a}
  \int_{\Omega}|u_n|^{\nu}|\phi| dx\rightarrow \int_{\Omega}|u|^{\nu}|\phi| dx,\,\,\,\text{as}\,\,\, n\rightarrow\infty
\end{equation}
and
\begin{equation}\label{strong3b}
  \int_{\Omega}|u_n|^{q}|\phi|(e^{\alpha u_n^2}-1)dx\rightarrow \int_{\Omega}|u|^{q}|\phi|(e^{\alpha u^2}-1)dx,\,\,\,\text{as}\,\,\, n\rightarrow\infty.
\end{equation}
Since $\nu>2-\frac{\mu}{4}>1$, we have $u_n\rightarrow u$ in $L^{2\nu}(\Omega)$, thus $|u_n|^\nu\rightarrow |u|^\nu$ in $L^{2}(\Omega)$. By the definition of weak convergence, we obtain \eqref{strong3a}. Furthermore, for fixed $q>0$, since  $t'=\frac{t}{t-1}$ large enough, we obtain $u_n\rightarrow u$ in $L^{qt'}(\Omega)$, thus $|u_n|^q\rightarrow |u|^q$ in $L^{t'}(\Omega)$, and
\begin{align*}
  &\Big|\int_{\Omega}|u_n|^{q}|\phi|(e^{\alpha u_n^2}-1)dx-\int_{\Omega}|u|^{q}|\phi|(e^{\alpha u^2}-1)dx\Big|\\
  &\leq\|\phi\|_\infty\int_{\Omega}\Big||u_n|^{q}-|u|^{q}\Big|(e^{\alpha u_n^2}-1)dx+\|\phi\|_\infty\int_{\Omega}|u|^{q}\Big|(e^{\alpha u_n^2}-1)-(e^{\alpha u^2}-1)\Big|dx\\
  &\leq \|\phi\|_\infty \Big(\int_{\Omega}\Big||u_n|^{q}-|u|^{q}\Big|^{t'}dx\Big)^{\frac{1}{t'}}\Big(\int_{\Omega}(e^{\alpha u_n^2}-1)^tdx\Big)^{\frac{1}{t}}\\
  &\quad+\|\phi\|_\infty\int_{\Omega}|u|^{q}\Big|(e^{\alpha u_n^2}-1)-(e^{\alpha u^2}-1)\Big|dx\rightarrow 0,\,\,\,\text{as}\,\,\, n\rightarrow\infty.
\end{align*}
Applying a variant of the Lebesgue dominated convergence theorem, we conclude the result.
\end{proof}

\begin{lemma}\label{Pohozaev}
If $u\in H^2(\mathbb{R}^4)$ is a critical point of $\mathcal J(u)$ on $S(c)$, then $u\in \mathcal{P}(c)$.
\end{lemma}
\begin{proof}
Consider a cut-off function $\varphi\in C_0^\infty(\mathbb{R}^4,[0,1])$ such that $\varphi(x)=1$ if $|x|\leq1$, $\varphi(x)=0$ if $|x|\geq2$.
For any fixed $\rho>0$, set $\widetilde{u}_\rho(x)=\varphi(\rho x)x\cdot\nabla u(x)$ as a test function of \eqref{e1.1} to obtain
\begin{equation*}
  \int_{\mathbb{R}^4}\Delta u \Delta\widetilde{u}_\rho dx+ \beta\int_{\mathbb{R}^4}\nabla u \cdot\nabla \widetilde{u}_\rho dx=\lambda \int_{\mathbb{R}^4} u\widetilde{u}_\rho dx+\int_{\mathbb{R}^4} (I_\mu*F(u))f(u)\widetilde{u}_\rho dx.
\end{equation*}
By \cite{MS2}, we know
\begin{equation*}
  \lim_{\rho\rightarrow0}\int_{\mathbb{R}^4}\nabla u\cdot \nabla \widetilde{u}_\rho dx=-\int_{\mathbb{R}^4}|\nabla u|^2 dx,\quad\lim_{\rho\rightarrow0}\int_{\mathbb{R}^4} u\widetilde{u}_\rho dx=-2\int_{\mathbb{R}^4}| u|^2 dx,
\end{equation*}
and
\begin{equation*}
 \lim_{\rho\rightarrow0} \int_{\mathbb{R}^4} (I_\mu*F(u))f(u)\widetilde{u}_\rho dx=-\frac{8-\mu}{2} \int_{\mathbb{R}^4} (I_\mu*F(u))F(u)dx.
\end{equation*}
Integrating by parts, we find that
\begin{align*}
  \int_{\mathbb{R}^4}\Delta u \Delta\widetilde{u}_\rho dx&=2\int_{\mathbb{R}^4}\varphi(\rho x)|\Delta u|^2 dx+\int_{\mathbb{R}^4}\varphi(\rho x)\Delta u(x\cdot \nabla (\Delta u))dx\\
  &=2\int_{\mathbb{R}^4}\varphi(\rho x)|\Delta u|^2 dx+\int_{\mathbb{R}^4}\varphi(\rho x)(x\cdot \nabla (\frac{|\Delta u|^2}{2}))dx
  \\&=2\int_{\mathbb{R}^4}\varphi(\rho x)|\Delta u|^2 dx-\int_{\mathbb{R}^4}\rho x\cdot\nabla \varphi(\rho x)\frac{|\Delta u|^2}{2}dx-\int_{\mathbb{R}^4}4\varphi(\rho x)\frac{|\Delta u|^2}{2}dx.
\end{align*}
The Lebesgue dominated convergence theorem implies that
\begin{equation*}
  \lim_{\rho\rightarrow0} \int_{\mathbb{R}^4}\Delta u \Delta\widetilde{u}_\rho dx=0.
\end{equation*}
Thus
\begin{equation}\label{po1}
  -\beta\int_{\mathbb{R}^4}|\nabla u|^2 dx=-2\lambda\int_{\mathbb{R}^4} |u|^2 dx-\frac{8-\mu}{2} \int_{\mathbb{R}^4} (I_\mu*F(u))F(u)dx.
\end{equation}
If $u\in H^2(\mathbb{R}^4)$ is a critical point of $\mathcal J(u)$ on $S(c)$, testing \eqref{e1.1} with $u$, we have
\begin{equation*}
  \int_{\mathbb{R}^4}|\Delta u |^2\ dx+\beta\int_{\mathbb{R}^4}|\nabla u|^2 dx=\lambda \int_{\mathbb{R}^4} |u|^2 dx+\int_{\mathbb{R}^4} (I_\mu*F(u))f(u)u dx.
\end{equation*}
This together with \eqref{po1} yields that
 \begin{align*}
  &2\int_{\mathbb{R}^4}|\Delta u |^2\ dx+\beta\int_{\mathbb{R}^4}|\nabla u|^2 dx+\frac{8-\mu}{2} \int_{\mathbb{R}^4} (I_\mu*F(u))F(u)dx-2\int_{\mathbb{R}^4} (I_\mu*F(u))f(u)u dx=0.
\end{align*}
Thus $u\in \mathcal{P}(c)$.
\end{proof}

To prove Theorem \ref{th1}, we consider the functional $H^2(\mathbb{R}^4)\times \mathbb{R}\rightarrow \mathbb{R}$ by
\begin{align*}
&\widetilde{\mathcal J}(u,s)=\mathcal J(\mathcal{H}(u,s))
=\frac{e^{4s}}{2}\int_{\mathbb{R}^4}|\Delta u|^2dx+\frac{\beta e^{2s}}{2}\int_{\mathbb{R}^4}|\nabla u|^2dx-\frac{e^{(\mu-8)s}}{2}\int_{\mathbb{R}^4}(I_\mu*F(e^{2s}u))F(e^{2s}u)dx
\end{align*}
Following  by \cite{Wil}, we recall that for any $c\in (0,\sqrt{\frac{4-\mu}{4}})$, the tangent space of $S(c)$ at $u$ is defined by
\begin{equation*}
  T_u=\Big\{v\in H^2(\mathbb{R}^4):\int_{\mathbb{R}^4}uvdx=0\Big\},
\end{equation*}
and the tangent space of $S(c)\times \mathbb{R}$ at $(u,s)$ is defined by
\begin{equation*}
  \widetilde{T}_{u,s}=\Big\{(v,t)\in H^2(\mathbb{R}^4)\times \mathbb{R}:\int_{\mathbb{R}^4}uvdx=0\Big\}.
\end{equation*}

\section{The minimax approach}\label{minimax}

In this section, we prove that $\widetilde{\mathcal J}$ on $S(c)\times \mathbb{R}$ possesses mountain-pass geometrical structure.

\begin{lemma}\label{minimax1}
Assume that $(f_1)$, $(f_2)$ and $(f_4)$ hold. If $c\in (0,\sqrt{\frac{4-\mu}{4}})$, then for any fixed $u\in S(c)$, we have

$(i)$ $\|\Delta \mathcal{H}(u,s)\|_2\rightarrow 0^+$ and $\mathcal{J}(\mathcal{H}(u,s))\rightarrow0^+$ as $s\rightarrow-\infty$;

$(ii)$ $\|\Delta \mathcal{H}(u,s)\|_2\rightarrow +\infty$ and $\mathcal{J}(\mathcal{H}(u,s))\rightarrow-\infty$ as $s\rightarrow+\infty$.
\end{lemma}
\begin{proof} A straightforward calculation shows that for any $q>2$,
\begin{equation*}
 \|\mathcal{H}(u,s)\|_2=c,\,\,\, \|\Delta \mathcal{H}(u,s)\|_2=e^{2s} \|\Delta u\|_2,\,\,\, \|\nabla \mathcal{H}(u,s)\|_2=e^{s} \|\nabla u\|_2,\,\,\,\|\mathcal{H}(u,s)\|_q=e^{\frac{2(q-2)s}{q}}\|u\|_q.
\end{equation*}
From the above equalities, we get $\|\Delta \mathcal{H}(u,s)\|_2\rightarrow 0$ as $s\rightarrow-\infty$, and $\|\Delta \mathcal{H}(u,s)\|_2\rightarrow +\infty$ as $s\rightarrow+\infty$. Furthermore, since $c\in (0,\sqrt{\frac{4-\mu}{4}})$, there exist $s_1<<0$ and $m\in (0,\frac{4-\mu}{4})$ such that
\begin{equation*}
  \|\mathcal{H}(u,s)\|^2= e^{4s} \|\Delta u\|_2^2+e^{2s} \|\nabla u\|_2^2
  +c^2\leq m<\frac{4-\mu}{4},\quad \text{for all $s\leq s_1$}.
\end{equation*}
By Proposition \ref{HLS}, fix $q>0$, for any $\xi>0$, $\alpha>32\pi^2$ close to $32\pi^2$ and $t>1$ close to $1$ with  $\frac{8\alpha mt}{8-\mu}\leq 32\pi^2$, there exists a constant $C_\xi>0$ such that 
\begin{align}\label{close1}
 &\int_{\mathbb{R}^4}(I_\mu*F(\mathcal{H}(u,s)))F(\mathcal{H}(u,s))dx\leq C\|F(\mathcal{H}(u,s))\|_{\frac{8}{8-\mu}}^2\nonumber\\
\leq & C\bigg(\int_{\mathbb{R}^4}\Big(\xi|\mathcal{H}(u,s)|^{\nu+1}+C_\xi|\mathcal{H}(u,s)|^{q+1}(e^{\alpha \mathcal{H}^2(u,s)}-1)\Big)^{\frac{8}{8-\mu}}dx\bigg)^{\frac{8-\mu}{4}}\nonumber\\
\leq& C \bigg(\int_{\mathbb{R}^4}\Big(|\mathcal{H}(u,s)|^{\frac{8(\nu+1)}{8-\mu}}+|\mathcal{H}(u,s)|^{\frac{8(q+1)}{8-\mu}}(e^{\frac{8\alpha \mathcal{H}^2(u,s)}{8-\mu}}-1)\Big)dx\bigg)^{\frac{8-\mu}{4}}\nonumber\\
\leq& C\|\mathcal{H}(u,s)\|_{\frac{8(\nu+1)}{8-\mu}}^{2(\nu+1)}+C\bigg(\int_{\mathbb{R}^4}|\mathcal{H}(u,s)|^{\frac{8(q+1)}{8-\mu}}(e^{\frac{8\alpha \mathcal{H}^2(u,s)}{8-\mu}}-1)dx\bigg)^{\frac{8-\mu}{4}}\nonumber\\
\leq&  C\|\mathcal{H}(u,s)\|_{\frac{8(\nu+1)}{8-\mu}}^{2(\nu+1)}+C\Big(\int_{\mathbb{R}^4}|\mathcal{H}(u,s)|^{\frac{8(q+1)t'}{8-\mu}}dx\Big)^{\frac{8-\mu}{4t'}}\Big(\int_{\mathbb{R}^4}(e^{\frac{8\alpha mt}{8-\mu}(\frac{\mathcal{H}(u,s)}{\|\mathcal{H}(u,s)\|})^2}-1)dx\Big)^{\frac{8-\mu}{4t}}\nonumber\\
\leq& C\|\mathcal{H}(u,s)\|_{\frac{8(\nu+1)}{8-\mu}}^{2(\nu+1)}+C\|\mathcal{H}(u,s)\|_{\frac{8(q+1)t'}{8-\mu}}^{2(q+1)}\nonumber\\
=&Ce^{(4\nu+\mu-4)s}\|u\|_{\frac{8(\nu+1)}{8-\mu}}^{2(\nu+1)}+Ce^{(4q+4-\frac{8-\mu}{t'})s}\|u\|_{\frac{8(q+1)t'}{8-\mu}}^{2(q+1)},\quad \text{for all $s\leq s_1$}.
\end{align}
Since $\beta\geq0$, $\nu>2-\frac{\mu}{4}$, $q>0$, and $t'=\frac{t}{t-1}$ large enough, it follows from \eqref{close1} that
\begin{equation*}
  \mathcal{J}(\mathcal{H}(u,s))\geq \frac{e^{4s}}{2}\|\Delta u\|_2^2-Ce^{(4\nu+\mu-4)s}\|u\|_{\frac{8(\nu+1)}{8-\mu}}^{2(\nu+1)}-Ce^{(4q+4-\frac{8-\mu}{t'})s}\|u\|_{\frac{8(q+1)t'}{8-\mu}}^{2(q+1)}\rightarrow 0^+, \quad\text{as $s\rightarrow-\infty$.}
\end{equation*}
By $(f_4)$,
\begin{align*}
\mathcal J(\mathcal{H}(u,s))
&\leq\frac{1}{2}\int_{\mathbb{R}^4}|\Delta \mathcal{H}(u,s)|^2dx+\frac{\beta}{2}\int_{\mathbb{R}^4}|\nabla \mathcal{H}(u,s)|^2dx-\frac{\tau^2}{2}\int_{\mathbb{R}^4}(I_\mu*|\mathcal{H}(u,s)|^p)|\mathcal{H}(u,s)|^pdx\\
&=\frac{e^{4s}}{2}\int_{\mathbb{R}^4}|\Delta u|^2dx+\frac{\beta e^{2s}}{2}\int_{\mathbb{R}^4}|\nabla u|^2dx-\frac{\tau^2e^{(4p+\mu-8)s}}{2}\int_{\mathbb{R}^4}(I_\mu*|u|^p)|u|^pdx.
\end{align*}
Since $p>3-\frac{\mu}{4}$, the above inequality yields that $\mathcal{J}(\mathcal{H}(u,s))\rightarrow-\infty$ as $s\rightarrow+\infty$.
\end{proof}
\begin{lemma}\label{minimax2}
Assume that $(f_1)-(f_3)$ hold. Then for any $c\in (0,\sqrt{\frac{4-\mu}{4}})$, there exist $0<k_1<k_2$ such that
\begin{equation*}
  0<\inf\limits_{u\in \mathcal{A}} \mathcal J(u)\leq \sup\limits_{u\in \mathcal{A}} \mathcal J(u)<\inf\limits_{u\in \mathcal{B}} \mathcal J(u)
\end{equation*}
with
\begin{equation*}
  \mathcal{A}=\Big\{u\in S(c):\|\Delta u\|_2\leq k_1\Big\},\quad\mathcal{B}=\Big\{u\in S(c):\|\Delta u\|_2=k_2\Big\}.
\end{equation*}
\end{lemma}

\begin{proof}
For any $u\in S(c)$ with $\|\Delta u\|_2\leq\frac{1}{3}(\sqrt{\frac{4-\mu}{4}}-c)$, by \eqref{GNinequality} and \eqref{close1}, fix $q>0$, for any $t>1$ close to $1$,
\begin{align*}
  \int_{\mathbb{R}^4}(I_\mu*F(u))F(u)dx\leq C\|u\|_{\frac{8(\nu+1)}{8-\mu}}^{2(\nu+1)}+C\|u\|_{\frac{8(q+1)t'}{8-\mu}}^{2(q+1)}\leq Cc^{4-\frac{\mu}{2}}\|\Delta u\|_2^{2\nu+\frac{\mu}{2}-2}+C c^{\frac{8-\mu}{2t'}} \|\Delta u\|_2^{2(q+1)-\frac{8-\mu}{2t'}}.
\end{align*}
By $(f_3)$, $F(u)\geq0$ for any $u\in H^2(\mathbb{R}^4)$. Since $\beta\geq0$, $\nu>2-\frac{\mu}{4}$, $q>0$, and $t'=\frac{t}{t-1}$ large enough, it follows that there exist $0<k_2\leq \frac{1}{3}(\sqrt{\frac{4-\mu}{4}}-c)$ small enough and $\rho>0$ such that
\begin{equation*}
  \mathcal{J}(v)\geq \frac{1}{2}\|\Delta v\|_2^2-Cc^{4-\frac{\mu}{2}}\|\Delta v\|_2^{2\nu+\frac{\mu}{2}-2}-C c^{\frac{8-\mu}{2t'}} \|\Delta v\|_2^{2(q+1)-\frac{8-\mu}{2t'}}\geq\rho>0,\quad\text{for all $v\in \mathcal{B}$}.
\end{equation*}
On the other hand, by \eqref{ineq}, we have
\begin{equation*}
  \mathcal{J}(u)\leq \frac{1}{2}\|\Delta u\|_2^2+\frac{\beta c}{2}\|\Delta u\|_2+Cc^{4-\frac{\mu}{2}}\|\Delta u\|_2^{2\nu+\frac{\mu}{2}-2}+C c^{\frac{8-\mu}{2t'}} \|\Delta u\|_2^{2(q+1)-\frac{8-\mu}{2t'}},
\end{equation*}
which implies that $\sup\limits_{u\in \mathcal{A}}\mathcal{J}(u)<\rho$ for any $k_1\in(0,k_2)$ small enough. And there exists $\sigma>0$ such that
\begin{equation*}
  \mathcal{J}(u)\geq \frac{1}{2}\|\Delta u\|_2^2-Cc^{4-\frac{\mu}{2}}\|\Delta u\|_2^{2\nu+\frac{\mu}{2}-2}-C c^{\frac{8-\mu}{2t'}} \|\Delta u\|_2^{2(q+1)-\frac{8-\mu}{2t'}}\geq\sigma>0,\quad\text{for all $u\in \mathcal{A}$}.
\end{equation*}
\end{proof}

\begin{lemma}\label{minimax3}
Let $k_1,k_2$ be defined in Lemma \ref{minimax2}. Then there exist $\hat{u},\tilde{u}\in S(c)$ such that

$(i)$ $\|\Delta \hat{u}\|_2\leq k_1$;

$(ii)$ $\|\Delta \tilde{u}\|_2> k_2$;

$(iii)$ $\mathcal J(\hat{u})>0>\mathcal J(\tilde{u})$.
\\Moreover, setting
\begin{equation*}
  \widetilde{m}_\tau(c)=\inf\limits_{\widetilde{h}\in \widetilde{\Gamma}_c}\max\limits _{t\in[0,1]}\widetilde{\mathcal J}({\widetilde{h}(t)})
\end{equation*}
with
\begin{equation*}
  \widetilde{\Gamma}_c=\Big\{\widetilde{h}\in C([0,1],S(c)\times\mathbb{R} ):\widetilde{h}(0)=(\hat{u},0),\widetilde{h}(1)=(\tilde{u},0)\Big\},
\end{equation*}
and
\begin{equation*}
  m_\tau(c)=\inf\limits_{h\in \Gamma_c}\max\limits _{t\in[0,1]}\mathcal J({h(t)})
\end{equation*}
with
\begin{equation*}
  \Gamma_c=\Big\{h\in C([0,1],S(c)):h(0)=\hat{u},h(1)=\tilde{u}\Big\},
\end{equation*}
then we have
\begin{equation*}
   \widetilde{m}_\tau(c)=m_\tau(c)\geq\max\{\mathcal J(\hat{u}),\mathcal J(\tilde{u})\}>0.
\end{equation*}
\end{lemma}

\begin{proof}
For any fixed $u_0\in S(c)$, by Lemmas \ref{minimax1} and \ref{minimax2}, there exist two numbers $s_1<<-1$ and $s_2>>1$ such $\hat{u}=\mathcal{H}(u_0,s_1)$ and $\tilde{u}=\mathcal{H}(u_0,s_2)$ satisfy $(i)-(iii)$. For any $\widetilde{h}\in \widetilde{\Gamma}_c$, we write it into
\begin{equation*}
  \widetilde{h}(t)=(\widetilde{h}_1(t),\widetilde{h}_2(t))\in S(c)\times\mathbb{R}.
\end{equation*}
Setting $h(t)=\mathcal{H}(\widetilde{h}_1(t),\widetilde{h}_2(t))$, then $h(t)\in \Gamma_c$ and
\begin{equation*}
  \max\limits _{t\in[0,1]}\widetilde{\mathcal J}({\widetilde{h}(t)})=\max\limits _{t\in[0,1]}\mathcal J({h(t)})\geq  m_\tau(c).
\end{equation*}
By the arbitrariness of $\widetilde{h}\in \widetilde{\Gamma}_c$, we get $\widetilde{m}_\tau(c)\geq m_\tau(c)$.

On the other hand, for any $h\in \Gamma_c$, if we set $\widetilde{h}(t)=(h(t),0)$, then $\widetilde{h}(t)\in \widetilde{\Gamma}_c$ and
\begin{equation*}
  \widetilde{m}_\tau(c)\leq \max\limits _{t\in[0,1]}\widetilde{\mathcal J}({\widetilde{h}(t)})=\max\limits _{t\in[0,1]}\mathcal J({h(t)}).
\end{equation*}
By the arbitrariness of $h\in {\Gamma}_c$, we get $\widetilde{m}_\tau(c)\leq m_\tau(c)$. Hence, we have $\widetilde{m}_\tau(c)= m_\tau(c)$, and $m_\tau(c)\geq\max\{\mathcal J(\hat{u}),\mathcal J(\tilde{u})\}$ follows from the definition of $m_\tau(c)$.
\end{proof}

Learning from \cite[Proposition 2.2]{Jeanjean}, by a standard Ekeland's variational principle and pseudo-gradient flow, we have the following proposition, which gives the existence of the $(PS)_{\widetilde{m}_\tau(c)}$ sequence for $\widetilde{\mathcal J}(u,s)$ on $S(c)\times \mathbb{R}$.
\begin{proposition}\label{psse}
Let $\widetilde{h}_n\subset \widetilde{\Gamma}_c$ be such that
\begin{equation*}
  \max\limits _{t\in[0,1]}\widetilde{\mathcal J}({\widetilde{h}_n(t)})\leq \widetilde{m}_\tau(c)+\frac{1}{n}.
\end{equation*}
Then there exists a sequence $\{(v_n,s_n)\}\subset S(c)\times\mathbb{R} $ such that as $n\rightarrow\infty$,

$(i)$ $\widetilde{\mathcal J}(v_n,s_n)\rightarrow \widetilde{m}_\tau(c)$;

$(ii)$ $\widetilde{\mathcal J}'|_{S(c)\times\mathbb{R}}(v_n,s_n)\rightarrow0$, i.e.,
\begin{equation*}
  \partial_s\widetilde{\mathcal J}(v_n,s_n)\rightarrow0 \quad \text{and}\quad\langle\partial_v\widetilde{\mathcal J}(v_n,s_n),\widetilde{\varphi}\rangle\rightarrow0
\end{equation*}
for all
\begin{equation*}
  \widetilde{\varphi}\in T_{v_n}=\Big\{\widetilde{\varphi}=(\widetilde{\varphi}_1,\widetilde{\varphi}_2)\in H^2(\mathbb{R}^4)\times \mathbb{R}:\int_{\mathbb{R}^4}v_n\widetilde{\varphi}_1\,dx=0\Big\}.
\end{equation*}
\end{proposition}
\begin{lemma}\label{pssequence}
For the sequence $\{(v_n,s_n)\}\subset S(c)\times\mathbb{R}$ obtained in Proposition \ref{psse}, setting $u_n=\mathcal{H}(v_n,s_n)$, then as $n\rightarrow\infty$, we have

$(i)$ ${\mathcal J}(u_n)\rightarrow {m}_\tau(c)$;

$(ii)$ $P(u_n)\rightarrow0$;

$(iii)$ ${\mathcal J}'|_{S(c)}(u_n)\rightarrow0$, i.e.,
\begin{equation*}
  \langle{\mathcal J}'(u_n),{\varphi}\rangle\rightarrow0
\qquad
\mbox{for\ all}\ \ \ \
  {\varphi}\in T_{u_n}=\Big\{{\varphi}\in H^2(\mathbb{R}^4):\int_{\mathbb{R}^4}u_n{\varphi}\,dx=0\Big\}.
\end{equation*}
\end{lemma}

\begin{proof}
For $(i)$, since ${\mathcal J}(u_n)=\widetilde{\mathcal J}(v_n,s_n)$ and ${m}_\tau(c)= \widetilde{m}_\tau(c)$, we get the conclusion.

For $(ii)$, first, we have
\begin{align*}
  \partial_s\widetilde{\mathcal J}(v_n,s_n)=&\partial_s\bigg[\frac{ e^{4s_n}}{2}\int_{\mathbb{R}^4}|\Delta v_n|^2dx+\frac{\beta e^{2s_n}}{2}\int_{\mathbb{R}^4}|\nabla v_n|^2dx-\frac{e^{(\mu-8)s_n}}{2}\int_{\mathbb{R}^4}(I_\mu*F(e^{2s_n}v_n))F(e^{2s_n}v_n)dx\bigg]\\
=&2 e^{4s_n}\int_{\mathbb{R}^4}|\Delta v_n|^2dx+\beta e^{2s_n}\int_{\mathbb{R}^4}|\nabla v_n|^2dx
+\frac{8-\mu}{2}e^{(\mu-8)s_n}\int_{\mathbb{R}^4}(I_\mu*F(e^{2s_n}v_n))F(e^{2s_n}v_n)dx
\\&-2e^{(\mu-8)s_n}\int_{\mathbb{R}^4}(I_\mu*F(e^{2s_n}v_n))f(e^{2s_n}v_n)e^{2s_n}v_ndx=P(u_n).
\end{align*}
Thus $P(u_n)\rightarrow0$ as $n\rightarrow\infty$.

For $(iii)$, for any $\widetilde{\varphi}=(\widetilde{\varphi}_1,\widetilde{\varphi}_2)\in T_{v_n}$,
\begin{align*}
\langle\partial_v\widetilde{\mathcal J}(v_n,s_n),\widetilde{\varphi}\rangle=&e^{4s_n}\int_{\mathbb{R}^4}\widetilde{\varphi} _1\Delta^2  v_n dx-\beta e^{2s_n}\int_{\mathbb{R}^4}\widetilde{\varphi}_1\Delta v_n dx-e^{(\mu-8)s_n}\int_{\mathbb{R}^4}(I_\mu*F(e^{2s_n}v_n))f(e^{2s_n}v_n)e^{2s_n}\widetilde{\varphi}_1dx.
  \end{align*}
On the other hand,  \begin{align*}
                       \langle{\mathcal J}'(u_n),{\varphi}\rangle=&\int_{\mathbb{R}^4}\varphi \Delta ^2u_n dx-\beta\int_{\mathbb{R}^4}\varphi \Delta u_n dx-\int_{\mathbb{R}^4}(I_\mu*F(u_n))f(u_n)\varphi dx\\
=& e^{2s_n}\int_{\mathbb{R}^4}\varphi(e^{-s_n}x)\Delta ^2v_n  dx-\beta\int_{\mathbb{R}^4}\varphi (e^{-s_n}x)\Delta v_ndx\\
&-e^{(\mu-8)s_n}\int_{\mathbb{R}^4}(I_\mu*F(e^{2s_n}v_n))f(e^{2s_n}v_n))\varphi(e^{-s_n}x) dx.
                    \end{align*}
Taking $\varphi(e^{-s_n}x) =e^{2s_n}\widetilde{\varphi}_1$, then $\langle{\mathcal J}'(u_n),{\varphi}\rangle\rightarrow 0$ as $n\rightarrow\infty$, and $\varphi(x) =e^{2s_n}\widetilde{\varphi}_1(e^{s_n}x)$. If we can prove $\varphi \in T_{u_n}$, we get $(iii)$. In fact, it follows from the following equality:
\begin{equation*}
  \int_{\mathbb{R}^4}u_n{\varphi}dx=\int_{\mathbb{R}^4}e^{2s_n}v_n(e^{s_n}x)e^{2s_n}\widetilde{\varphi}_1(e^{s_n}x)dx= \int_{\mathbb{R}^4}v_n\widetilde{{\varphi}}_1dx=0.
\end{equation*}
\end{proof}
\begin{lemma}\label{small}
Assume that $(f_1)-(f_4)$ hold, then $\lim\limits_{\tau\rightarrow\infty}m_\tau(c)=0$.
\end{lemma}
\begin{proof}
For any fixed $u_0\in H^2(\mathbb{R}^4)$, as the proof of Lemma \ref{minimax3},  for any $t\in [0,1]$, $h_0(t)=\mathcal{H}(u_0,(1-t)s_1+ts_2)$ is a path in $\Gamma_c$. Hence, by $(f_4)$,
\begin{align*}
  m_\tau(c)\leq \max\limits_{t\in[0,1]}\mathcal J({h_0(t)})\leq& \max\limits_{\kappa>0}\Big\{\frac{\kappa^4}{2}\|\Delta u_0\|_2^2+\frac{\beta \kappa^2}{2}\|\nabla u_0\|_2^2
-\frac{\tau^2\kappa^{(4p+\mu-8)}}{2}\int_{\mathbb{R}^4}(I_\mu*|u_0|^p)|u_0|^pdx\Big\}\\
\leq&\max\limits_{\kappa>0}\Big\{\frac{\kappa^4}{2}\|\Delta u_0\|_2^2-\frac{\tau^2\kappa^{(4p+\mu-8)}}{4}\int_{\mathbb{R}^4}(I_\mu*|u_0|^p)|u_0|^pdx\Big\}
\\&+\max\limits_{\kappa>0}\Big\{\frac{\beta \kappa^2}{2}\|\nabla u_0\|_2^2
-\frac{\tau^2\kappa^{(4p+\mu-8)}}{4}\int_{\mathbb{R}^4}(I_\mu*|u_0|^p)|u_0|^pdx\Big\}\\=&C\Big(\frac{1}{\tau^2}\Big)^{\frac{4}{4p+\mu-12}}+C\Big(\frac{1}{\tau^2}\Big)^{\frac{2}{4p+\mu-10}}.
\end{align*}
This together with $p>3-\frac{\mu}{4}$ yields $\lim\limits_{\tau\rightarrow\infty}m_\tau(c)=0$.
\end{proof}
\begin{lemma}\label{lowerenergy}
Assume that $(f_1)-(f_4)$ hold and $c\in (0,\sqrt{\frac{4-\mu}{4}})$. For the sequence $\{u_n\}$ obtained in Lemma \ref{pssequence}, 
there exists $\tau_1>0$ such that for any $\tau\geq\tau_1$,
\begin{equation*}
  \limsup\limits_{n\rightarrow\infty}\|\Delta u_n\|_2<\sqrt{\frac{4-\mu}{4}}-c.
\end{equation*}
\end{lemma}
\begin{proof}
By ${\mathcal J}(u_n)\rightarrow {m}_\tau(c)$, $P(u_n)\rightarrow0$ as $n\rightarrow\infty$ and $(f_3)$, we get
\begin{align*}
  {\mathcal J}(u_n)-\frac{1}{4}P(u_n)&=\frac{\beta }{4}\int_{\mathbb{R}^4}|\nabla u_n|^2dx+\frac{\mu-12}{8}\int_{\mathbb{R}^4}(I_\mu*F(u_n))F(u_n)dx+\frac{1}{2}\int_{\mathbb{R}^4}(I_\mu*F(u_n))f(u_n)u_n dx\nonumber\\
  &\geq\frac{\beta }{4}\int_{\mathbb{R}^4}|\nabla u_n|^2dx+\frac{4\theta+\mu-12}{8}\int_{\mathbb{R}^4}(I_\mu*F(u_n))F(u_n)dx.
\end{align*}
Since $\beta\geq0$, by $\theta>3-\frac{\mu}{4}$, we have
\begin{equation*}
 \limsup\limits_{n\rightarrow\infty} \int_{\mathbb{R}^4} (I_\mu*F(u_n))F(u_n) dx\leq \frac{8 m_\tau(c)}{4\theta+\mu-12}.
\end{equation*}
And immediately, we obtain
\begin{equation*}
  \limsup\limits_{n\rightarrow\infty}\frac{1}{2}\int_{\mathbb{R}^4}|\Delta u_n |^2 dx\leq \limsup\limits_{n\rightarrow\infty}\frac{1}{2}\int_{\mathbb{R}^4} (I_\mu*F(u_n))F(u_n) dx+m_\tau(c)\leq \frac{4\theta+\mu-8 }{4\theta+\mu-12}m_\tau(c).
\end{equation*}
Combining this with Lemma \ref{small}, we derive the conclusion.
\end{proof}

For the sequence $\{u_n\}$ obtained in Lemma \ref{pssequence}, by Lemma \ref{lowerenergy}, $\{u_n\}$ is bounded in $H^2(\mathbb{R}^4)$, up to a subsequence, we assume that $u_n\rightharpoonup u_c$ in $H^2(\mathbb{R}^4)$.
 Furthermore, by ${\mathcal J}'|_{S(c)}(u_n)\rightarrow0$ as $n\rightarrow\infty$ and Lagrange multiplier rule, there exists $\{\lambda_n\}\subset \mathbb{R}$ such that
\begin{equation}\label{lagrange}
   \Delta^2u_n-\beta\Delta u_n=\lambda_n u_n+(I_\mu*F(u_n))f(u_n)+o_n(1).
\end{equation}
\begin{lemma}\label{up}
Assume that $(f_1)-(f_4)$ hold and $c\in (0,\sqrt{\frac{4-\mu}{4}})$. Up to a subsequence and up to translations in $\mathbb{R}^4$, $u_n\rightharpoonup u_c\neq0$ in $H^2(\mathbb{R}^4)$ for any $\tau\geq\tau_1$.
\end{lemma}
\begin{proof}
If $\tau\geq \tau_1$, by Lemma \ref{lowerenergy}, up to a subsequence, we assume that $\sup\limits_{n\in \mathbb{N}}\|u_n\|^2<m<\frac{4-\mu}{4}$. Since $\{u_n\}$ is bounded in $H^2(\mathbb{R}^4)$,
for any $r>0$, let
\begin{equation*}
  \varrho:=\limsup\limits_{n\rightarrow\infty}\Big(\sup\limits_{y\in \mathbb{R}^4}\int_{B(y,r)}|u_n|^2dx\Big).
\end{equation*}

If $\varrho>0$, then there exists $\{y_n\}\subset \mathbb{R}^4$  such that $\int_{B(y_n,1)}|u_n|^2dx>\frac{\varrho}{2}$, i.e., $\int_{B(0,1)}|u_n(x-y_n)|^2dx>\frac{\varrho}{2}$. Up to a subsequence and up to translations in $\mathbb{R}^4$, $u_n\rightharpoonup u_c\neq0$ in $H^2(\mathbb{R}^4)$.

If $\varrho=0$, by \cite[Lemma 1.21]{Wil}, $u_n\rightarrow0$ in $L^p(\mathbb{R}^4)$ for any $p>2$. As a consequence,
arguing as \eqref{close1}, by $\nu>2-\frac{\mu}{4}$, $q>0$ and $t'=\frac{t}{t-1}$ large enough, we have
 \begin{align*}
  \int_{\mathbb{R}^4}(I_\mu*F(u_n))f(u_n)u_ndx&\leq C\|u_n\|_{\frac{8(\nu+1)}{8-\mu}}^{2(\nu+1)}+C\|u_n\|_{\frac{8(q+1)t'}{8-\mu}}^{2(q+1)}
 \rightarrow0,\quad\text{as $n\rightarrow\infty$.}
\end{align*}
From the above equality and $P(u_n)\rightarrow0$, we deduce that $\|\Delta u_n\|_2\rightarrow 0$ as $n\rightarrow\infty$, hence $(f_3)$ implies $\lim\limits_{n\rightarrow\infty}{\mathcal J}(u_n)\leq  0$, which is an absurd, since $m_\tau(c)>0$.
\end{proof}

\begin{lemma}\label{lambda}
Assume that $(f_1)-(f_4)$ hold and $c\in (0,\sqrt{\frac{4-\mu}{4}})$. Then $\{\lambda_n\}$ is bounded and
\begin{equation*}
\limsup\limits_{n\rightarrow\infty}\lambda_n=\limsup\limits_{n\rightarrow\infty}\frac{\beta}{2c^2}\int_{\mathbb{R}^4}|\nabla u_n|^2dx- \liminf\limits_{n\rightarrow\infty}\frac{8-\mu}{4c^2}\int_{\mathbb{R}^4}(I_\mu*F(u_n))F(u_n)dx.
\end{equation*}
\end{lemma}
\begin{proof}
Testing \eqref{lagrange} with $u_n$, we have
\begin{equation*}
   \int_{\mathbb{R}^4}|\Delta u_n|^2dx+\beta\int_{\mathbb{R}^4}|\nabla u_n|^2dx=\lambda_n \int_{\mathbb{R}^4}|u_n|^2dx+\int_{\mathbb{R}^4}(I_\mu*F(u_n))f(u_n)u_ndx+o_n(1).
\end{equation*}
Combing this with $P(u_n)\rightarrow0$ as $n\rightarrow\infty$ lead to
\begin{equation*}
  \lambda_n c^2=\frac{\beta}{2}\int_{\mathbb{R}^4}|\nabla u_n|^2dx-\frac{8-\mu}{4}\int_{\mathbb{R}^4}(I_\mu*F(u_n))F(u_n)dx+o_n(1).
\end{equation*}
Thus,\begin{equation*}
   \limsup\limits_{n\rightarrow\infty}\lambda_n = \limsup\limits_{n\rightarrow\infty}\frac{\beta}{2c^2}\int_{\mathbb{R}^4}|\nabla u_n|^2dx- \liminf\limits_{n\rightarrow\infty}\frac{8-\mu}{4c^2}\int_{\mathbb{R}^4}(I_\mu*F(u_n))F(u_n)dx.
\end{equation*}
From (\ref{ineq}), we find
\begin{equation*}
  \limsup\limits_{n\rightarrow\infty}|\lambda_n| \leq \limsup\limits_{n\rightarrow\infty}\Big(\frac{\beta}{2c}\Big(\int_{\mathbb{R}^4}|\Delta u_n|^2dx\Big)^{\frac{1}{2}}+\frac{8-\mu}{4c^2}\int_{\mathbb{R}^4}(I_\mu*F(u_n))F(u_n)dx\Big).
\end{equation*}
This together with Lemma \ref{lowerenergy} yields that
\begin{equation*}
  \limsup\limits_{n\rightarrow\infty}|\lambda_n| \leq \frac{\beta}{2c}\sqrt{\frac{2(4\theta+\mu-8 )}{(4\theta+\mu-12)}m_\tau(c)}+\frac{2(8-\mu) }{(4\theta+\mu-12)c^2}m_\tau(c).
\end{equation*}
\end{proof}
Since $u_c\neq0$, by Lemmas \ref{small}, \ref{lowerenergy}, \ref{lambda}, and Fatou lemma, 
we find that there exists $\tau_2>0$ such that for any $\tau\geq\tau_2$,
\begin{align*}
  \limsup\limits_{n\rightarrow\infty}\frac{\beta}{2c^2}\int_{\mathbb{R}^4}|\nabla u_n|^2dx& \leq\frac{\beta}{2c}\sqrt{\frac{2(4\theta+\mu-8 )}{(4\theta+\mu-12)}m_\tau(c)}\\&< \frac{8-\mu}{4c^2}\int_{\mathbb{R}^4}(I_\mu*F(u_c))F(u_c)dx\\&\leq\liminf\limits_{n\rightarrow\infty}\frac{8-\mu}{4c^2}\int_{\mathbb{R}^4}(I_\mu*F(u_n))F(u_n)dx,
\end{align*}
which implies that $\limsup\limits_{n\rightarrow\infty}\lambda_n<0$ for any $\tau\geq\tau_2$. Since $\{\lambda_n\}$ is bounded, up to a subsequence, we assume that $\lambda_n\rightarrow\lambda_c<0$ as $n\rightarrow\infty$ for any $\tau\geq\tau_2$.

\begin{lemma}\label{equi}
Assume that $f$ satisfies $(f_1)-(f_4)$ and $c\in(0,\sqrt{\frac{4-\mu}{4}})$. If $(f_5)$ or $(f_6)$ holds,
then for any fixed $u\in S(c)$,
the function ${\mathcal J}(\mathcal{H}(u,s))$ achieves its maximum with positive level at a unique point $s_u\in \mathbb{R}$ such that $\mathcal{H}(u,s_u) \in \mathcal{P}(c)$.
\end{lemma}
\begin{proof}
By Lemma \ref{minimax1}, we have
\begin{equation*}
  \lim\limits_{s\rightarrow-\infty}\mathcal{J}(\mathcal{H}(u,s))=0^+\quad\text{and}\lim\limits_{s\rightarrow+\infty}\mathcal{J}(\mathcal{H}(u,s))=-\infty.
\end{equation*}
Therefore, there exists $s_u\in \mathbb{R}$ such that $P(\mathcal{H}(u,s_u))=\frac{d}{ds}\mathcal{J}(\mathcal{H}(u,s))|_{s=s_u}=0$, and $\mathcal{J}(\mathcal{H}(u,s_u))>0$.
In the following, we prove the uniqueness of $s_u$ for any $u\in S(c)$.

Case 1: If $(f_1)-(f_5)$ hold. Taking into account that $\frac{d}{ds}\mathcal{J}(\mathcal{H}(u,s))|_{s=s_u}=0$, using $\beta\geq0$ and $(f_5)$, we deduce that
\begin{align*}
  \frac{d^2}{ds^2}\mathcal{J}(\mathcal{H}(u,s))|_{s=s_u}=&
8e^{4s_u}\int_{\mathbb{R}^4}|\Delta u|^2dx+2\beta e^{2s_u}\int_{\mathbb{R}^4}|\nabla u|^2dx
\\&-\frac{(8-\mu)^2}{2}e^{(\mu-8)s_u}\int_{\mathbb{R}^4}(I_\mu*F(e^{2s_u}u))F(e^{2s_u}u)dx
\\&+(28-4\mu)e^{(\mu-8)s_u}\int_{\mathbb{R}^4}(I_\mu*F(e^{2s_u}u))f(e^{2s_u}u)e^{2s_u}udx\\
&-4e^{(\mu-8)s_u}\int_{\mathbb{R}^4}(I_\mu*f(e^{2s_u}u)e^{2s_u}u)f(e^{2s_u}u)e^{2s_u}udx\\&
-4e^{(\mu-8)s_u}\int_{\mathbb{R}^4}(I_\mu*F(e^{2s_u}u))f'(e^{2s_u}u)e^{4s_u}u^2dx\\
=&-2\beta \int_{\mathbb{R}^4}|\nabla \mathcal{H}(u,s_u)|^2dx-(8-\mu)(6-\frac{\mu}{2})\int_{\mathbb{R}^4}(I_\mu*F(\mathcal{H}(u,s_u)))F(\mathcal{H}(u,s_u))dx\\
&+(36-4\mu)\int_{\mathbb{R}^4}(I_\mu*F(\mathcal{H}(u,s_u)))f(\mathcal{H}(u,s_u))\mathcal{H}(u,s_u)dx\\
&-4\int_{\mathbb{R}^4}(I_\mu*f(\mathcal{H}(u,s_u))\mathcal{H}(u,s_u))f(\mathcal{H}(u,s_u))\mathcal{H}(u,s_u)dx\\&
-4\int_{\mathbb{R}^4}(I_\mu*F(\mathcal{H}(u,s_u)))f'(\mathcal{H}(u,s_u))\mathcal{H}^2(u,s_u)dx \\
=&-2\beta \int_{\mathbb{R}^4}|\nabla \mathcal{H}(u,s_u)|^2dx +4\int_{\mathbb{R}^4}\int_{\mathbb{R}^4}\frac{A}{|x-y|^\mu}dxdy<0,
\end{align*}
where
\begin{align*}
  A=&(3-\frac{\mu}{4})F(\mathcal{H}(u(y),s_{u(y)}))\overline{F}(\mathcal{H}(u(x),s_{u(x)}))-F(\mathcal{H}(u(y),s_{u(y)}))\overline{F}'(\mathcal{H}(u(x),s_{u(x)}))\mathcal{H}(u(x),s_{u(x)})\\&
  -\overline{F}(\mathcal{H}(u(y),s_{u(y)}))(\overline{F}(\mathcal{H}(u(x),s_{u(x)}))-F(\mathcal{H}(u(x),s_{u(x)})))<0,
\end{align*}
this prove the uniqueness of $s_u$.

Case 2: If $(f_1)-(f_4)$ and $(f_6)$ hold. Since
\begin{align*}
P(\mathcal{H}(u,s))=&e^{4s}\Big[2 \int_{\mathbb{R}^4}|\Delta u|^2dx+ \frac{\beta}{e^{2s}}\int_{\mathbb{R}^4}|\nabla u|^2dx -2\int_{\mathbb{R}^4}\Big(I_\mu*\frac{F(e^{2s}u)}{(e^{2s})^{3-\frac{\mu}{4}}}\Big)\frac{\overline{F}(e^{2s}u)}{(e^{2s})^{3-\frac{\mu}{4}}}\Big]dx.
\end{align*}
Denote
\begin{equation*}
\psi(s)=\int_{\mathbb{R}^4}\Big(I_\mu*\frac{F(e^{2s}u)}{(e^{2s})^{3-\frac{\mu}{4}}}\Big)\frac{\overline{F}(e^{2s}u)}{(e^{2s})^{3-\frac{\mu}{4}}}dx.
\end{equation*}
For any $t\in \mathbb{R}\backslash \{0\}$, from $(f_3)$ and $(f_6)$, we see that $\frac{F(st)}{s^{3-\frac{\mu}{4}}}$ is increasing in $s\in (0,+\infty)$ and $\frac{\overline{F}(st)}{s^{3-\frac{\mu}{4}}}$ is non-decreasing in $s\in (0,+\infty)$. This fact implies $\psi(s)$ is increasing in $s\in \mathbb{R}$ and there is at most one $s_u\in \mathbb{R}$ such that $\mathcal{H}(u,s_u) \in \mathcal{P}(c)$. This ends the proof.
\end{proof}

\begin{remark}\label{control}
{\rm
 By Lemma \ref{equi}, we claim that $m_\tau(c)=E(c)$. In fact, by Lemma \ref{equi}, we know $\mathcal{P}(c)\neq \emptyset$. For any fixed $u_0\in \mathcal{P}(c)$, as in Lemma \ref{small}, we define a path $h_0(t)=\mathcal{H}(u_0,(1-t)s_1+ts_2)$  in $\Gamma_c$ for any $t\in [0,1]$. Then by Lemma \ref{equi}, $\max\limits_{t\in[0,1]}\mathcal{J}(h_0(t))=\mathcal{J}(u_0)$, which implies that
\begin{equation*}
  E(c)\geq \inf\limits_{h\in \Gamma _c}\max\limits_{t\in[0,1]}\mathcal{J}(h(t))=m_\tau(c).
\end{equation*}
On the other hand, since $S(c)\backslash \mathcal{P}(c)$ has two components given by $R^{+}:=\{u\in S(c):P(u)>0\}$ and $R^{-}:=\{u\in S(c):P(u)<0\}$. Since $\|\Delta \hat{u}\|_2\leq k_1$ and $k_1>0$ small enough, with a similar proof of Lemma \ref{minimax2}, we deduce that $\hat{u}\in R^{+}$. By $\mathcal{J}(\tilde{u})<0$ and $(f_3)$, we see that
\begin{align*}
  P(\tilde{u})&<(6-\frac{\mu}{2})\int_{\mathbb{R}^4}(I_\mu*F(\tilde{u})F(\tilde{u})dx-2\int_{\mathbb{R}^4}(I_\mu*F(\tilde{u})
  f(\tilde{u})\tilde{u}dx\\
  &\leq (6-\frac{\mu}{2}-2\theta)\int_{\mathbb{R}^4}(I_\mu*F(\tilde{u})F(\tilde{u})dx<0,
\end{align*}
which shows that $\tilde{u}\in R^{-}$. Hence, any path in $\Gamma_c$ must intersect $\mathcal{P}(c)$, and we have
\begin{equation*}
  E(c)\leq  \inf\limits_{h\in \Gamma _c}\max\limits_{t\in[0,1]}\mathcal{J}(h(t))=m_\tau(c).
\end{equation*}
}
\end{remark}

\section{Proof of the main result}\label{main}

\noindent
{\it Proof of Theorem \ref{th1}:}
Since $f$ is odd in $\mathbb{R}$, it's easy to find that $\int_{\mathbb{R}^4}(I_\mu*F(u))F(u)dx$ and $\int_{\mathbb{R}^4}(I_\mu*F(u))f(u)udx$ are even in $H^2(\mathbb{R}^4)$. Then we can see that all the above conclusions can be repeated word by word in $H^2_{rad}(\mathbb{R}^4)$, and without loss of generality, we assume that $u_n\geq0$ for all $n\in \mathbb{N}^+$.

From \eqref{lagrange}, $\lambda_n\rightarrow \lambda_c$ as $n\rightarrow\infty$, we know
\begin{equation*}
  \int_{\mathbb{R}^4}|\Delta u_n|^2dx+\beta\int_{\mathbb{R}^4}|\nabla u_n|^2dx=\lambda_c \int_{\mathbb{R}^4}|u_n|^2dx+\int_{\mathbb{R}^4}(I_\mu*F(u_n))f(u_n)u_ndx+o_n(1).
\end{equation*}
By Lemma \ref{strong3}, we can see that $u_c$ is a weak solution of \eqref{e1.6}, i.e.,
\begin{equation*}
  \Delta^2u_c-\beta\Delta u_c=\lambda_c u_c+(I_\mu*F(u_c))f(u_c),
\end{equation*}
thus
\begin{equation*}
   \int_{\mathbb{R}^4}|\Delta u_c|^2dx+\beta\int_{\mathbb{R}^4}|\nabla u_c|^2dx=\lambda_c \int_{\mathbb{R}^4}|u_c|^2dx+\int_{\mathbb{R}^4}(I_\mu*F(u_c))f(u_c)u_cdx.
\end{equation*}
Choosing $\tau\geq\tau^*=\max\{\tau_1,\tau_2\}$, the above equalities together with Lemma \ref{strong2} and $\lambda_c<0$ gives
\begin{equation*}
  \lim\limits_{n\rightarrow\infty}\int_{\mathbb{R}^4}|\Delta u_n|^2dx=\int_{\mathbb{R}^4}|\Delta u_c|^2dx, \quad\lim\limits_{n\rightarrow\infty}\int_{\mathbb{R}^4}|u_n|^2dx=\int_{\mathbb{R}^4}|u_c|^2dx
\end{equation*}
which implies $u_n\rightarrow u_c$ in $H^2_{rad}(\mathbb{R}^4)$. Therefore, $u_c$ is a non-negative radial ground state solution of \eqref{e1.6}.

\subsection*{Acknowledgments}
The authors have been supported by National Natural Science Foundation of China 11971392 and Natural Science Foundation of Chongqing, China cstc2021ycjh-bgzxm0115.

\end{document}